\documentclass[a4paper,11pt,twoside,reqno]{amsart}
\usepackage[utf8]{inputenc}
\usepackage{amsmath,amsthm,amssymb}

\usepackage[margin=3.2cm]{geometry}

\usepackage{listings}
\usepackage{xcolor}

\lstdefinestyle{mypython}{
  language=Python,
  basicstyle=\ttfamily\small,
  keywordstyle=\color{blue},
  commentstyle=\color{teal!70!black},
  stringstyle=\color{red!60!black},
  showstringspaces=false,
  keepspaces=true,
  columns=fullflexible,
  breaklines=true,
  breakatwhitespace=false,
  frame=single,
  numbers=left,
  numberstyle=\scriptsize\color{gray},
  stepnumber=1,
  numbersep=8pt,
  tabsize=4,
  captionpos=b
}

\usepackage[pagebackref=true]{hyperref}
\usepackage[alphabetic,initials]{amsrefs}

\def\ab{\mathbf{a}}
\def\bb{\mathbf{b}}
\def\cb{\mathbf{c}}
\def\Ab{\mathbf{A}}
\def\Bb{\mathbf{B}}
\def\Cb{\mathbf{C}}
\newcommand\ddd{\,\mathrm{d}}

\newcommand\RR{\mathbb{R}}
\newcommand\NN{\mathbb{N}}
\newcommand\ZZ{\mathbb{Z}}

\renewcommand{\leq}{\leqslant}
\renewcommand{\le}{\leqslant}
\renewcommand{\geq}{\geqslant}
\renewcommand{\ge}{\geqslant}
\newcommand\eps{\varepsilon}
\newcommand\ve{\varepsilon}

\theoremstyle{plain}
\newtheorem*{theorem*}{Theorem}
\newtheorem{prop}{Proposition}
\newtheorem{lemma}[prop]{Lemma}
\newtheorem{theorem}[prop]{Theorem}
\theoremstyle{remark}

\makeatletter
\@namedef{subjclassname@2020}{\textup{2020} Mathematics Subject Classification}
\makeatother
\DeclareMathOperator\rad{rad}

\numberwithin{prop}{section}
\numberwithin{equation}{section}

\title{Bounds on the exceptional set in the $abc$ conjecture}

\author{C. Bernert}
\author{T. Browning}
\address{Institute for Science and Technology Austria\\
Am Campus 1\\
3400 Klosterneuburg\\
Austria}
\email{christian.bernert@ist.ac.at, tdb@ist.ac.at}

\author{J. D. Lichtman}
\address{
Department of Mathematics\\ 
Stanford University\\ 
450 Jane Stanford Way\\
Stanford, CA 94305-2125\\
USA}
\email{jdlichtman@stanford.edu}

\author{J. Ter\"av\"ainen}
\address{Department of Pure Mathematics and Mathematical Statistics\\ University of Cambridge\\
Cambridge CB3 0WB\\ UK}
\email{joni.p.teravainen@gmail.com}

\date{\today}

\begin{document}

\begin{abstract}
We study solutions to the equation $a+b=c$, where $a,b,c$ form a triple of coprime natural numbers.
The $abc$ conjecture asserts that, for any $\ve>0$, such triples satisfy $\rad(abc) \ge c^{1-\ve}$ with finitely many exceptions. 
In this article we obtain a power-saving bound on the size of the exceptional set of triples. 
The proof is based on a combination of upper bounds for the density of integer points on certain high-dimensional varieties,  
coming from the geometry of numbers and from Fourier analysis.
\end{abstract}

\maketitle

\setcounter{tocdepth}{1}
\tableofcontents

\section{Introduction}

In this paper, we use tools from analytic number theory to estimate the number of triples of a given height satisfying the $abc$ conjecture.  
Associated to any non-zero  integer $n$ is its radical 
\begin{align*}
\rad(n)=\prod_{p\mid n}p.    
\end{align*}
We say that a triple $(a,b,c)\in \NN^3$
with $\gcd(a,b,c)=1$ is 
 an $abc$ \emph{triple of exponent} $\lambda$ if
\begin{align*}
a+b=c,\quad \rad(abc) < c^{\lambda}.    
\end{align*}
The well-known $abc$ conjecture of Masser and Oesterl\'e asserts that, for any $\lambda<1$, there are only finitely many $abc$ triples of exponent $\lambda$. See also the work of Robert, Stewart and Tenenbaum~\cite{tenenbaum} for a refinement of the $abc$ conjecture. 
The   best unconditional result is  due to Stewart and Yu~\cite{styu}, who  have shown that only finitely many $abc$ triples satisfy $\rad(abc) < (\log c)^{3-\eps}$. Recently, Pasten~\cite{hector} has proved a new subexponential bound, assuming that  $a<c^{1-\eps}$, via a connection to Shimura curves. 
In this paper we shall focus on counting the number $N_\lambda(X)$ of $abc$ triples of exponent $\lambda$ in a box $[1,X]^3$, as $X\to \infty$. 

Given $\lambda>0$, 
an old result of
de Bruijn~\cite{debruijn} implies that 
\begin{equation}\label{eq:bruijn}
\#\left\{n\leq x : \rad(n)\leq x^\lambda\right\}=O_\ve(x^{\lambda+\ve}),
\end{equation}
for any $\ve>0$. Any triple $(a,b,c)$ counted by 
$N_\lambda(X)$ must satisfy $\rad(abc)< X^{\lambda}$, and so 
we must have
$\min\{\rad(a)\rad(b),\rad(b)\rad(c),\rad(c)\rad(a)\}< X^{2\lambda/3}$, since $a,b,c$ are pairwise coprime.
An application of \eqref{eq:bruijn} now leads to the following ``trivial bound''.

\begin{prop}\label{prop:trivial}
Let $\lambda>0$. Then 
$
N_\lambda(X)=O_\ve( X^{2\lambda/3+\ve}),
$
for any $\ve>0$.
\end{prop}

Note that for $\lambda\geq 3$ the counting function $N_\lambda(X)$ is asymptotic to a constant multiple of $X^2$, as $X\to \infty$,  and so it is only worth  studying when  $\lambda <3$. 
The primary goal of this paper is to give the first power-saving improvement over the simple bound
in Proposition \ref{prop:trivial}
 for values of $\lambda\in (0,2)$. In 2000, it was asked by Mazur \cite{mazur} whether 
$N_\lambda(X)$ has exact order $X^{\lambda-1}$, for a  fixed $\lambda > 1$. 
This was answered positively by Kane~\cite{kane} when 
  $\lambda \ge 2$,  who used  the determinant method  to prove that 
$$
X^{\lambda-1+\varepsilon}\ll 
N_{\lambda}(X) \ll_{\varepsilon} X^{\lambda-1+\varepsilon},
$$
for any $\ve>0$.
(In fact, the lower bound is valid for any $\lambda> 1$ and is originally due to unpublished work of Granville.) We assume henceforth that $\lambda\in (0,2]$. The following is our main result.  

\begin{theorem} \label{thm:main_general}
Let  $\lambda \in (0,2]$. For any $\ve>0$, we have
    \[N_{\lambda}(X) \ll_{\varepsilon} X^{\frac{23\lambda+3}{40}+\varepsilon}.\]
    In particular, we have 
    $
    N_\lambda(X) \ll_{\lambda} X^{0.65}
    $,
if $\lambda\in (0,1)$.
\end{theorem}

Our result is non-trivial for a relatively large range of $\lambda$, since it follows that 
$$
N_{\lambda}(X)=o(X^{2\lambda/3}),
$$
for all $\lambda \in (9/11,2]$. Given the relationship to the $abc$ conjecture, the case $\lambda\in (0,1)$ is of special interest. For such $\lambda$, 
a more refined combination of our bounds leads to a rather complicated optimisation problem, which will yield the following improvement. 

\begin{theorem} \label{thm:main_optimal}
Let $\lambda \in (0,1)$. Then $N_\lambda(X) \ll_{\ve,\lambda} X^{0.6+\ve}$, for any $\ve>0$.
\end{theorem}

This result is non-trivial for $\lambda\in (9/10,1)$.
The present paper merges an earlier article 
of the last three authors \cite{blt24}, in which additional bounds coming from 
the determinant method and the theory of Thue equations were used to produce a worse bound, together with  its subsequent improvement and simplification by the first author 
 \cite{bernert25}.
Since our bounds do not use the determinant method, it is plausible that further improvements are available, but we expect that a serious numerical improvement will be likely to require new ideas.

Besides our work, we are not aware of any general estimates for $N_\lambda(X)$ when $\lambda<1$, 
 beyond  Proposition \ref{prop:trivial}. Nonetheless, there do exist    specific Diophantine equations which are covered by the $abc$ conjecture and where bounds have been given for the number of solutions. 
For example, it follows from work of Darmon and 
Granville \cite{DG} that there are only finitely many coprime integer solutions to the Diophantine equation $x^p+y^q=z^r$, when $p,q,r\in \mathbb{N}$ are given and satisfy $1/p+1/q+1/r<1$.

\subsection*{Notation}
We shall use $x\sim X$ to denote $x\in (X,2X]$.
All statements involving a number $\varepsilon$ are meant to hold for all sufficiently small constants $\varepsilon>0$.
In particular, we will frequently use the ``divisor bound'' by which we mean the fact that the number of divisors of a positive integer $n$ is $O_{\varepsilon}(n^{\varepsilon})$.
For a real number $\alpha$ we write $e(\alpha)=e^{2\pi  i\alpha}$ for the corresponding complex number on the unit circle.

\subsection*{Acknowledgements}
We are grateful to Bhavik Mehta for his assistance with using Lean to verify parts of an earlier version of this article. 
While working on this paper
the third author was supported by an NSF Postdoctoral Fellowship.
The fourth author was supported by the European Union's Horizon
Europe research and innovation programme under Marie Sk\l{}odowska-Curie grant agreement no. 101058904, by ERC grant agreement no. 101162746, and by Academy of Finland grant no. 362303. 
 This material is based upon work supported by a grant from the Institute
for Advanced Study School of Mathematics.

\section{An anatomic reduction}

Let us denote by $N_{\lambda}^*(X)$ the number of tuples $(a,b,c)$ counted by $N_{\lambda}(X)$ satisfying additionally $c \ge X/2$ and $\rad(abc)\geq c^{\lambda}/2$.
By a dyadic decomposition and the pigeonhole principle, we clearly have
\[
N_{\lambda}(X) \ll (\log X)^2 \max_{1\leq Y\leq X}\max_{\mu\leq \lambda}N_{\mu}^*(Y).
\]
Thus it will suffice  to bound $N_{\lambda}^*(X)$.

The main goal of this section is to bound 
$N_{\lambda}^*(X)$ 
in terms of the  number of solutions to certain monomial Diophantine equations. In order to state it, we need to introduce  the quantity
$$
B_d(a_0,b_0,c_0, \Ab,\Bb,\Cb)
=  \#\left\{
(\ab,\bb,\cb) \in \NN^{3d}
\;:\;
\begin{array}{l}
a_i \sim A_i,\, b_i \sim B_i,\, c_i \sim C_i,\\
a_0 \prod_{i=1}^d a_i^i+b_0\prod_{i=1}^d b_i^i=c_0 \prod_{i=1}^d c_i^i,\\
\gcd\left(\prod_{i=0}^d a_i, \prod_{i=0}^d b_i, \prod_{i=0}^d c_i\right)=1
\end{array}
\right\},
$$
for $d \in \NN$, a vector $(a_0,b_0,c_0) \in \NN^3$ and parameters 
$$
\Ab=(A_1,\dots,A_d), \Bb=(B_1,\dots,B_d), \Cb=(C_1,\dots,C_d) \in \NN^d.
$$
Bearing this notation in mind, we shall prove the following bound.

\begin{prop}\label{prop:reduction}
    For any $\varepsilon_0>0$, there exists a positive integer $d=d(\varepsilon_0)$ such that
    \[
    N_{\lambda}^*(X) \ll_{\varepsilon_0} X^{\varepsilon_0} \cdot \max B_d(a_0,b_0,c_0,\Ab,\Bb,\Cb),
    \]
    where the maximum runs over pairwise coprime $1\le a_0,b_0,c_0 \le X^{\varepsilon_0}$ and dyadic parameters $A_i,B_i,C_i \in [1,X]$ with
    \[
    \frac{1}{2^{d+1}}X^{\lambda}\leq 
    \prod_{j=1}^d A_jB_jC_j \le X^{\lambda+3\varepsilon_0/4}\]
    and
    \[\prod_{j=1}^d A_j^j \le X, \quad \prod_{j=1}^d B_j^j \le X, \quad X^{1-\varepsilon_0/2} \le \prod_{j=1}^d C_j^j \le X.\]
\end{prop}

This anatomic decomposition of the $abc$ equation is based on the following lemma, which provides such a decomposition for a single variable and which can easily be used to recover the bound 
\eqref{eq:bruijn}.

\begin{lemma}\label{lem:factorization}
    For any $\varepsilon_0>0$, there exists a positive integer $d=d(\varepsilon_0)\geq 6$ such that any positive integer $n$ has a factorization
    \[n=n_0 \prod_{j=1}^d n_j^j\]
    with positive integers $n_0,n_1,\dots,n_d$ such that $n_1,\dots,n_d$ are pairwise coprime, and satisfying $n_0 \le n^{\varepsilon_0}$ and
    \[n^{-\varepsilon_0} \prod_{j=1}^d n_j \le \rad(n) \le n^{\varepsilon_0} \prod_{j=1}^d n_j.\]
\end{lemma}

\begin{proof}
    We may assume that $\varepsilon_0<1$ for otherwise the claim is trivial with $d=6$.
        Fix $d=6\lceil \varepsilon_0^{-2}\rceil$ and $k=\lceil \varepsilon_0^{-1}\rceil$. For a given $n$, we consider 
    \[
    m_j:=\prod_{p^j\| n} p,
    \]
as  naive candidates for $n_j$.
We need to modify this construction to be able to truncate our products at $j \le d$. To this end, for $1 \le j \le d$, we define
    \[
    n_j:=\begin{cases} m_j & j \ne k\\ m_j \prod_{i>d} m_i^{\lfloor i/k\rfloor} & j=k \end{cases}
    \]
    as well as
    \[
    n_0:=\prod_{j>d} m_j^{j-k\lfloor j/k\rfloor}.
    \]
    Clearly, we have
    \[n=\prod_{j} m_j^j=n_0 \prod_{j=1}^d n_j^j
    \]
    and, since the $m_i$ are pairwise coprime, so are $n_1,\dots,n_d$. Moreover,
    \[
    n_0 \le \prod_{j>d} m_j^k \le n^{k/d} \le n^{\varepsilon_0},
    \]
    since $m_j^k\leq (m_j^{j/d})^k=(m_j^{j})^{k/d}$.
    Moving on to the radical, we have
    \[
    \rad(n) \le \rad(n_0)\prod_{j=1}^d \rad(n_j) \le n^{\varepsilon_0} \prod_{j=1}^d n_j.
    \]
On the other hand, in view of 
    \[
    n_k=m_k \prod_{j>d} m_j^{\lfloor j/k\rfloor} \le \left(m_k^k \prod_{j>d} m_j^j\right)^{1/k} \le n^{1/k} \le n^{\varepsilon_0},
    \]
    we find that
    \[
    \rad(n)=\prod_{j \ge 1} m_j \ge \prod_{j \le d, j \ne k} n_j \ge n^{-\varepsilon_0} \prod_{j=1}^d n_j,
    \]
    as required.
\end{proof}
\begin{proof}[Proof of Proposition~\ref{prop:reduction}]
    We may assume that $X$ is sufficiently large, for otherwise the claim is trivially true. We will choose $d=d(\varepsilon_0)=6\lceil (\varepsilon_0/4)^{-2}\rceil$.
    
    Suppose that $(a,b,c)$ is a triple counted by $N_{\lambda}^*(X)$. By Lemma~\ref{lem:factorization} we obtain factorizations
    \[
    a=a_0 \prod_{j=1}^d a_j^j, \quad b=b_0 \prod_{j=1}^d b_j^j, \quad c=c_0 \prod_{j=1}^d c_j^j,
    \]
    with $a_0,b_0,c_0 \le X^{\varepsilon_0/4}$. By assumption, $a+b=c$ and $\gcd(a,b,c)=1$, so that $a,b,c$ are pairwise coprime and hence also $\prod_{i=0}^d a_i$, $\prod_{i=0}^d b_i$ and $\prod_{i=0}^d c_i$ are pairwise coprime.

    Choosing $A_i,B_i,C_i$ to be the unique powers of two such that $a_i \sim A_i, b_i \sim B_i, c_i \sim C_i$, we find that
    \[
\frac{1}{2^{d+1}}X^\lambda \leq     \frac{\rad(abc)}{2^d}
\leq     \prod_{j=1}^d A_jB_jC_j \le \prod_{j=1}^d a_jb_jc_j \le \rad(abc)X^{3\varepsilon_0/4}  \le X^{\lambda+3\varepsilon_0/4}
    \]
and
    \[\prod_{j=1}^d A_j^j \le \prod_{j=1}^d a_j^j \le X,\]
together with similar upper bounds for $\prod_{j=1}^d B_j^j$ and $\prod_{j=1}^d C_j^j$. 
For the lower bound, we note that 
    \[\prod_{j=1}^d C_j^j \ge 2^{-d} \prod_{j=1}^d c_j^j \ge X^{1-\varepsilon_0/2},\]
    provided that $X$ is sufficiently large.

    This shows that $(a,b,c)$ is counted by $B_d(a_0,b_0,c_0,\Ab,\Bb,\Cb)$ for a suitable  choice of $A_i,B_i,C_i$ as powers of two. Since there
    are $O(X^{3\ve_0/4})$ choices for $a_0,b_0,c_0$ and at most $(\log X)^{3d} \ll_{\varepsilon_0} X^{\varepsilon_0/4}$ 
    choices for $A_i,B_i,C_i$,  the claim  follows from the pigeonhole principle.    
\end{proof}

In the following two sections, we will study bounds for $B_d(a_0,b_0,c_0,\Ab,\Bb,\Cb)$ for given parameters $a_0,b_0,c_0,A_i,B_i,C_i$ satisfying the conditions in Proposition~\ref{prop:reduction}, for a given value of $\varepsilon_0$ and the corresponding value of $d=d(\varepsilon_0)$, as well as for a  sufficiently large parameter $X$.
To simplify notation, let us write 
$$
A=\prod_{i=1}^d A_i, \quad B=\prod_{i=1}^d B_i, \quad  C=\prod_{i=1}^d C_i.
$$

Note that by fixing two of the three sets of variables in the defining equation of 
$B_d(a_0,b_0,c_0,\Ab,\Bb,\Cb)$,
the remaining ones are determined up to a divisor function. We may thus record our first estimate, which 
recovers Proposition~\ref{prop:trivial}, since  $ABC \le X^{\lambda+3\varepsilon_0}$.

\begin{prop}\label{prop:triv2}
    We have
    \[B_d(a_0,b_0,c_0,\Ab,\Bb,\Cb) \ll_{\varepsilon} X^{\varepsilon}\min(AB,AC,BC).\]
\end{prop}

\section{Fourier analysis}

In this section, we bound $B_d(a_0,b_0,c_0,\Ab,\Bb,\Cb)$ using Fourier analysis. (Note that our  bound also recovers 
Proposition~\ref{prop:trivial}.)

\begin{prop}\label{prop:fourierbound}
    For $2 \le j,k,\ell \le d$, we have
    \[B_d(a_0,b_0,c_0,\Ab,\Bb,\Cb) \ll_{\varepsilon_0} \frac{X^{2\lambda/3+3\varepsilon_0}}{\left(\prod_{j \mid i} A_i \prod_{k \mid i} B_i \prod_{\ell \mid i} C_i\right)^{1/6}}.\]
\end{prop}

\begin{proof}
    By orthogonality of characters and Hölder's inequality, we have
    \[
    B_d(a_0,b_0,c_0,\Ab,\Bb,\Cb) \le \int_0^1 |S_1(\alpha)S_2(\alpha)S_3(-\alpha)| \ddd\alpha \le \left(\prod_{i=1}^3 \int_0^1 |S_i(\alpha)|^3 \ddd\alpha\right)^{1/3},
    \]
    with the exponential sums
    \[S_1(\alpha)=\sum_{a_i \sim A_i} e(\alpha a_0a_1a_2^2\cdots a_d^d), \quad S_2(\alpha)=\sum_{b_i \sim B_i} e(\alpha b_0b_1b_2^2\cdots b_d^d)\]
    and
    \[S_3(\alpha)=\sum_{c_i \sim C_i} e(\alpha c_0c_1c_2^2\cdots c_d^d).\]
    We claim that 
    \begin{equation}\label{eq:second}
        \int_0^1 |S_1(\alpha)|^2 \ddd\alpha \ll_{d,\varepsilon} X^{\varepsilon} A
    \end{equation}
    and 
    \begin{equation}\label{eq:fourth}
        \int_0^1 |S_1(\alpha)|^4 \ddd\alpha \ll_{d,\varepsilon} X^{\varepsilon} \frac{A^3}{\prod_{j \mid i} A_i}.
    \end{equation}
    Combined with an application of Cauchy-Schwarz and the analogous inequalities for $S_2$ and $S_3$, this will suffice to deduce the result.

    For the proof of the second moment bound \eqref{eq:second}, we note that, by orthogonality, the second moment counts the number of solutions to $a_1a_2^2\cdots a_d^d=a_1'a_2'^2\cdots a_d'^d$     with $a_i,a_i' \sim A_i$. But on fixing all the $a_i$, the $a_i'$ are determined up to $O(X^{\varepsilon})$ many choices by the divisor bound. This proves \eqref{eq:second}.
    
    Turning to the fourth moment bound \eqref{eq:fourth}, we begin by applying Cauchy-Schwarz to the sum over $a_i$ for $j \nmid i$ in $S_1(\alpha)$ to obtain
    \[|S_1(\alpha)|^2 \le \prod_{j \nmid i} A_i \sum_{a_i \sim A_i} \sum_{a_i' \sim A_i, j \mid i}  e\left(\alpha a_0 \left(\prod_{j \mid i} a_i^i-\prod_{j \mid i} (a_i')^i\right)\prod_{j \nmid i} a_i^i\right).\]
    Inserting this bound into our fourth moment and using orthogonality of characters, we find that the fourth moment is bounded by $\prod_{j \nmid i} A_i$ times the number of solutions to
    \[
    \left(\prod_{j \mid i} a_i^i- \prod_{j \mid i} a_i'^i\right)\prod_{j \nmid i} a_i^i=\prod_i (a_i'')^i-\prod_i (a_i''')^i
    \]
    with $a_i,a_i'', a_i''' \sim A_i$ and, for $j \mid i$, also $a_i' \sim A_i$.
    It therefore suffices to show that the number of such solutions is bounded by $O(X^{\varepsilon} A^2)$.
    This is certainly true for the cases where both sides are equal to zero, since we can choose the $a_i$ and $a_i''$ freely and then the $a_i'$ and $a_i'''$ are determined up to a divisor function.
    
    In the remaining cases, we have at most $A^2$ many choices for the $a_i''$ and $a_i'''$.  But then the left-hand side is fixed. In particular, the $a_i$ for $j \nmid i$ are determined up to a divisor function. 
 We now recall the well-known fact that, for a fixed non-zero integer $n\in \ZZ$, 
   the Diophantine equation $x^j-y^j=n$ has $O_\ve(|nX|^\ve)$ integer solutions in the box $[1,X]^2$. 
   (This can be proved using the divisor bound.)
   Thus it follows that 
    $a_i$ and $a_i'$ are determined up to a divisor function. This proves \eqref{eq:fourth} and therefore the proposition.
\end{proof}

\section{Geometry of numbers}

In this section, we bound $B_d(a_0,b_0,c_0,\Ab,\Bb,\Cb)$ by fixing some of the variables and treating the remaining equation as a linear equation  via the geometry of numbers. This leads to the following bound.

\begin{prop}\label{prop:geometrybound}
    For any sets $I,I',I'' \subset \{1,2,\dots,d\}$ we have
    \begin{align*}
        B_d(a_0,b_0,c_0,\Ab,\Bb,\Cb) \ll_{\varepsilon_0}~& \frac{X^{\lambda+4\varepsilon_0}}{\prod_{i \in I} A_i \prod_{i \in I'} B_i \prod_{i \in I''} C_i} \\
        &+X^{\lambda-1+6\varepsilon_0} \prod_{i \in I} A_i^{i-1} \prod_{i \in I'} B_i^{i-1} \prod_{i \in I''} C_i^{i-1}.\end{align*}
\end{prop}

\begin{proof}
    Let $(a_1,\dots,a_d,b_1,\dots,b_d, c_1,\dots,c_d)$ be a tuple counted by $B_d(a_0,b_0,c_0,\Ab,\Bb,\Cb)$. Define
    \begin{alignat*}{3}
    &a'=a_0\prod_{i\not\in I}a_i^i, \quad 
    &&b'=b_0\prod_{i\not\in I'}b_i^i, \quad 
    &&c'=c_0\prod_{i\not\in I''}c_i^i,\\
    &x=\prod_{i\in I}a_i^i, \qquad
    &&y=\prod_{i\in I'}b_i^i, \qquad 
    &&z=\prod_{i\in I''}c_i^i,\\
    &X=\prod_{i\in I}A_i^i, \qquad
    &&Y=\prod_{i\in I'}B_i^i, \qquad 
    &&Z=\prod_{i\in I''}C_i^i.
\end{alignat*}
Then $\gcd(a',b',c')=1$ and $\gcd(x,y,z)=1$ and we have
\[a'x+b'y=c'z.\]
By a classical bound from the geometry of numbers, as proved by Heath-Brown \cite{hb84}*{Lemma~3}, for example, for fixed $a',b',c'$, the number of primitive solutions $(x,y,z)$ is bounded by
\[\ll 1+\frac{XYZ}{\max\{|a'|X,|b'|Y,|c'|Z\}}.\]
Any such solution fixes $a_i$ for $i \in I$, $b_i$ for $i \in I'$ and $c_i$ for $i \in I''$ up to $O(X^{\varepsilon})$ many choices. Summing over the choices of the remaining variables, 
the $O(1)$-term contributes 
\begin{align*}
\ll_\ve X^\ve \prod_{i\not\in I}A_i\prod_{i\not\in I'}B_i\prod_{i\not\in I''}C_i
&=
X^\ve \frac{ABC}{\prod_{i\in I}A_i\prod_{i\in I'}B_i\prod_{i\in I''}C_i}\\
&\ll_{\ve_0} \frac{X^{\lambda+4\ve_0}}{\prod_{i\in I}A_i\prod_{i\in I'}B_i\prod_{i\in I''}C_i},
\end{align*}
since $ABC\leq X^{\lambda+3\ve_0}$, 
which is satisfactory. Similarly, the remaining term contributes
\begin{align*}
&\ll_\ve X^\ve \frac{\prod_{i\in I}A_i^i\prod_{i\in I'}B_i^i\prod_{i\in I''}C_i^i}{X^{1-2\ve}} \cdot
\prod_{i\not\in I}A_i\prod_{i\not\in I'}B_i\prod_{i\not\in I''}C_i\\
&\ll_{\ve_0} X^{-1+3\ve_0}ABC \prod_{i\in I}A_i^{i-1}\prod_{i\in I'}B_i^{i-1}\prod_{i\in I''}C_i^{i-1},
\end{align*}
which  is also satisfactory.
\end{proof}

\section{The simplest non-trivial bound}

In this section we prove Theorem~\ref{thm:main_general}. For convenience let us write $P_i=A_iB_iC_i$.

\begin{proof}[Proof of Theorem \ref{thm:main_general}]
Let $\lambda\in (0,2]$ and suppose that the theorem is false.
By Proposition~\ref{prop:reduction}, we can then find $\delta_0<\frac{11\lambda-9}{120}$ such that for any sufficiently small $\varepsilon_0>0$,  and for the corresponding value of $d$, there exist $a_0,b_0,c_0,\Ab,\Bb,\Cb$ with
    \[B_d(a_0,b_0,c_0,\Ab,\Bb,\Cb) \ge X^{2\lambda/3-\delta_0},\]
    for arbitrarily large $X$. Fix $\delta \in (\delta_0,\frac{11\lambda-9}{120})$.
        For sufficiently small $\varepsilon_0>0$, it follows from 
        Proposition~\ref{prop:geometrybound} with $I=I'=I''=\{1\}$ that we must have $P_1 \ll X^{\lambda/3+\delta}$. By Proposition~\ref{prop:fourierbound} with $j=k=\ell$, we must have $P_j \ll X^{6\delta}$ for all $j \ge 2$. However,
    \begin{equation*}
        X^{5\lambda-3} \ll \prod_i P_i^{5-i} \ll P_1^4P_2^3P_3^2P_4 \ll X^{4\lambda/3+40\delta},
    \end{equation*}
    which contradicts our upper bound for $\delta$.
\end{proof}

\section{A stronger result for $\lambda<1$}

In this section we describe how to prove Theorem~\ref{thm:main_optimal}. To begin with,
the following result is obtained by solving a mixed-integer linear program; the corresponding \texttt{Python} code is given in Section~\ref{s:code}.

\begin{prop}\label{prop:linear_program}
    Let $a,a_1,\dots,a_6,b,b_1,\dots,b_6,c,c_1,\dots,c_6,D \in \RR_{\geq 0}$. Let  
    $$
    p_i=a_i+b_i+c_i, \quad L=a+b+c. 
     $$
 Then $D \le 0.6$ if the  following list of $22$ inequalities hold:
    \begin{equation}\label{eq:1}
            L \le 1,
    \end{equation}

\medskip

    \begin{equation}\label{eq:2}
\left\{
\begin{aligned}
a_1+a_2+a_3+a_4+a_5+a_6 &\le a,\\
b_1+b_2+b_3+b_4+b_5+b_6 &\le b,\\
c_1+c_2+c_3+c_4+c_5+c_6 &\le c,
\end{aligned}
\right.
\end{equation}

\medskip

\begin{equation}\label{eq:3}
\left\{
\begin{aligned}
        7a-1 &\le 6a_1+5a_2+4a_3+3a_4+2a_5+a_6,\\
        7b-1 &\le 6b_1+5b_2+4b_3+3b_4+2b_5+b_6,\\
        7c-1 &\le 6c_1+5c_2+4c_3+3c_4+2c_5+c_6,
\end{aligned}
\right.
\end{equation}

\medskip

        \begin{equation}\label{eq:4}
        D \le \min(a+b,a+c,b+c),
    \end{equation}

\medskip

\begin{equation}\label{eq:5}
       \min(p_1+p_2+p_4,1-p_2-3p_4) \le L-D,
    \end{equation}

\medskip

\begin{equation}\label{eq:6}
\left\{
\begin{aligned}
        \min(p_1+p_2+a_3,1-p_2-2a_3) &\le L-D,\\
        \min(p_1+p_2+b_3,1-p_2-2b_3) &\le L-D,\\
        \min(p_1+p_2+c_3,1-p_2-2c_3) &\le L-D,
\end{aligned}
\right.
\end{equation}

\medskip

\begin{equation}\label{eq:7}
\left\{
\begin{aligned}
        \min(p_1+p_2+a_3+b_3,1-p_2-2a_3-2b_3) &\le L-D,\\
        \min(p_1+p_2+b_3+c_3,1-p_2-2b_3-2c_3) &\le L-D,\\
        \min(p_1+p_2+c_3+a_3,1-p_2-2c_3-2a_3) &\le L-D,
\end{aligned}
\right.
\end{equation}

\medskip

\begin{equation}\label{eq:8}
\left\{
\begin{aligned}
        \min(p_1+a_3+b_2+c_2,1-2a_3-b_2-c_2) &\le L-D,\\
        \min(p_1+a_2+b_3+c_2,1-a_2-2b_3-c_2) &\le L-D,\\
        \min(p_1+a_2+b_2+c_3,1-a_2-b_2-2c_3) &\le L-D,\\
\end{aligned}
\right.
\end{equation}

\medskip

\begin{equation}\label{eq:9}
\left\{
\begin{aligned}
        \min(p_1+a_2+b_3+c_3,1-a_2-2b_3-2c_3) &\le L-D,\\
        \min(p_1+a_3+b_2+c_3,1-2a_3-b_2-2c_3) &\le L-D,\\
        \min(p_1+a_3+b_3+c_2,1-2a_3-2b_3-c_2) &\le L-D,\\
\end{aligned}
\right.
\end{equation}

\medskip

\begin{equation}\label{eq:10}
\begin{split}
        \max(a_2+a_4+a_6,a_3+a_6,a_5)&+\max(b_2+b_4+b_6,b_3+b_6,b_5)\\
        &+\max(c_2+c_4+c_6,c_3+c_6,c_5) 
        \le 4L-6D.
\end{split}
\end{equation}
\end{prop}

\begin{proof}[Proof of Theorem~\ref{thm:main_optimal}]
We fix a choice of  $\lambda\in (0,1)$. 
Suppose that 
Theorem~\ref{thm:main_optimal} is false. Then, by Proposition~\ref{prop:reduction}, we can find $\delta_0>0.6$ such that for any $\varepsilon_0$ and the corresponding $d$, there are $a_0,b_0,c_0, \Ab,\Bb,\Cb$ with
    \[B_d(a_0,b_0,c_0,\Ab,\Bb,\Cb) \ge X^{\delta_0}\]
    for arbitrarily large $X$. Fix $D \in (0.6,\delta_0)$.

    Write $A=X^a, B=X^b, C=X^c$ as well as $A_i=X^{a_i}, B_i=X^{b_i}, C_i=X^{c_i}$. Then, as we will see below, for sufficiently small $\varepsilon_0>0$, the variables $a_i,b_i,c_i,a,b,c,D$ satisfy the hypotheses of Proposition~\ref{prop:linear_program} contradicting the fact that $D>0.6$.

Assuming that $\ve_0$ is sufficiently small, it follows that $ABC\leq X^{\lambda+3\ve_0}\leq X$. This establishes  \eqref{eq:1}. The 
inequalities in
  \eqref{eq:2} 
follow immediately from the definitions. 
The inequalities in \eqref{eq:3} also follow from the definitions and the observation that 
$$
\sum_{1\leq i\leq 6} ia_i \leq 1, \quad 
\sum_{1\leq i\leq 6} ib_i \leq 1, \quad
\sum_{1\leq i\leq 6} ic_i \leq 1,
$$
so that 
$$
7a-1\le 7a-\sum_{1\leq i\leq 6} ia_i\leq 6a_1+5a_2+4a_3+3a_4+2a_5+a_6.
$$
Similarly for the remaining two inequalities in \eqref{eq:3}.
The  inequality \eqref{eq:4} is just the trivial bound from Proposition~\ref{prop:triv2}.

    The next set of  inequalities originate from the geometry bound in Proposition~\ref{prop:geometrybound} with various choices of $I$, $I'$, $I''$. 
       More precisely, \eqref{eq:5} follows on taking
    \[
    I=I'=I''=\{1,2,4\};
    \]
the inequalities \eqref{eq:6} correspond to taking 
    \[
    I=I'=\{1,2\}, I''=\{1,2,3\}
    \]
    and its permutations;
the inequalities \eqref{eq:7} correspond to taking 
    \[
    I=\{1,2\}, I'=I''=\{1,2,3\}
    \]
    and its permutations; 
    the inequalities \eqref{eq:8} correspond to taking 
    \[
    I=I'=\{1,2\}, I''=\{1,3\}
    \]
    and its permutations; and finally
the inequalities \eqref{eq:9} correspond to taking 
    \[
    I=\{1,2\}, I'=I''=\{1,3\}
    \]
    and its permutations.
The final  inequality \eqref{eq:10} is the Fourier bound from Proposition~\ref{prop:fourierbound}.
\end{proof}

\appendix
\section{Python code for the linear program}\label{s:code}

\begin{lstlisting}[style=mypython]
import pulp

# define problem
prob = pulp.LpProblem("LP_with_OR_constraints", pulp.LpMaximize)

# Variables
p1 = pulp.LpVariable("p1", lowBound=0)
a1 = pulp.LpVariable("a1", lowBound=0)
b1 = pulp.LpVariable("b1", lowBound=0)
c1 = pulp.LpVariable("c1", lowBound=0)
p2 = pulp.LpVariable("p2", lowBound=0)
a2 = pulp.LpVariable("a2", lowBound=0)
b2 = pulp.LpVariable("b2", lowBound=0)
c2 = pulp.LpVariable("c2", lowBound=0)
p3 = pulp.LpVariable("p3", lowBound=0)
a3 = pulp.LpVariable("a3", lowBound=0)
b3 = pulp.LpVariable("b3", lowBound=0)
c3 = pulp.LpVariable("c3", lowBound=0)
p4 = pulp.LpVariable("p4", lowBound=0)
a4 = pulp.LpVariable("a4", lowBound=0)
b4 = pulp.LpVariable("b4", lowBound=0)
c4 = pulp.LpVariable("c4", lowBound=0)
p5 = pulp.LpVariable("p5", lowBound=0)
a5 = pulp.LpVariable("a5", lowBound=0)
b5 = pulp.LpVariable("b5", lowBound=0)
c5 = pulp.LpVariable("c5", lowBound=0)
p6 = pulp.LpVariable("p6", lowBound=0)
a6 = pulp.LpVariable("a6", lowBound=0)
b6 = pulp.LpVariable("b6", lowBound=0)
c6 = pulp.LpVariable("c6", lowBound=0)
a = pulp.LpVariable("a", lowBound=0)
b = pulp.LpVariable("b", lowBound=0)
c = pulp.LpVariable("c", lowBound=0)
A = pulp.LpVariable("A", lowBound=0)
B = pulp.LpVariable("B", lowBound=0)
C = pulp.LpVariable("C", lowBound=0)
D  = pulp.LpVariable("D")
L = pulp.LpVariable("L", lowBound=0)

# binary variables for OR conditions
z1 = pulp.LpVariable("z1", cat="Binary")
z2 = pulp.LpVariable("z2", cat="Binary")
z3 = pulp.LpVariable("z3", cat="Binary")
z4 = pulp.LpVariable("z4", cat="Binary")
z5 = pulp.LpVariable("z5", cat="Binary")
z6 = pulp.LpVariable("z6", cat="Binary")
z7 = pulp.LpVariable("z7", cat="Binary")
z8 = pulp.LpVariable("z8", cat="Binary")
z9 = pulp.LpVariable("z9", cat="Binary")
z10 = pulp.LpVariable("z10", cat="Binary")
z11 = pulp.LpVariable("z11", cat="Binary")
z12 = pulp.LpVariable("z12", cat="Binary")
z13 = pulp.LpVariable("z13", cat="Binary")
M = 100  # Big-M

# Goal function: Maximal exponent
prob += D

# Setup
prob += p1 == a1+b1+c1
prob += p2 == a2+b2+c2
prob += p3 == a3+b3+c3
prob += p4 == a4+b4+c4
prob += p5 == a5+b5+c5
prob += p6 == a6+b6+c6
prob += L == a + b + c

prob += L <= 1
prob += a1 + a2 + a3 + a4 + a5 + a6 <= a
prob += b1 + b2 + b3 + b4 + b5 + b6 <= b
prob += c1 + c2 + c3 + c4 + c5 + c6 <= c
prob += 7*a - 1 <= 6*a1 + 5*a2 + 4*a3 + 3*a4 + 2*a5 + a6
prob += 7*b - 1 <= 6*b1 + 5*b2 + 4*b3 + 3*b4 + 2*b5 + b6
prob += 7*c - 1 <= 6*c1 + 5*c2 + 4*c3 + 3*c4 + 2*c5 + c6

# Trivial bounds
prob += D <= b + c
prob += D <= a + c
prob += D <= a + b

# Geometry bounds

# 124, 124, 124
prob += p1 + p2 + p4 <= L - D + M*z1
prob += 1 - p2 - 3*p4 <= L - D + M*(1 - z1)

# 12, 12, 123
prob += p1 + p2 + a3 <= L - D + M*z2
prob += 1 - p2 - 2*a3 <= L - D + M*(1 - z2)
prob += p1 + p2 + b3 <= L - D + M*z3
prob += 1 - p2 - 2*b3 <= L - D + M*(1 - z3)
prob += p1 + p2 + c3 <= L - D + M*z4
prob += 1 - p2 - 2*c3 <= L - D + M*(1 - z4)

# 12, 123, 123
prob += p1 + p2 + a3 + b3 <= L - D + M*z5
prob += 1 - p2 - 2*a3 - 2*b3 <= L - D + M*(1 - z5)
prob += p1 + p2 + b3 + c3 <= L - D + M*z6
prob += 1 - p2 - 2*b3 - 2*c3 <= L - D + M*(1 - z6)
prob += p1 + p2 + c3 + a3 <= L - D + M*z7
prob += 1 - p2 - 2*c3 - 2*a3 <= L - D + M*(1 - z7)

# 12,12,13
prob += p1 + a2 + b3 + c2 <= L - D + M*z8
prob += 1 - a2 - 2*b3 - c2 <= L - D + M*(1 - z8)
prob += p1 + a2 + c3 + b2 <= L - D + M*z9
prob += 1 - a2 - 2*c3 - b2 <= L - D + M*(1 - z9)
prob += p1 + b2 + a3 + c2 <= L - D + M*z10
prob += 1 - b2 - 2*a3 - c2 <= L - D + M*(1 - z10)

# 12,13,13
prob += p1 + b2 + c3 + a3 <= L - D + M*z11
prob += 1 - b2 - 2*c3 - 2*a3 <= L - D + M*(1 - z11)
prob += p1 + c2 + a3 + b3 <= L - D + M*z12
prob += 1 - c2 - 2*a3 - 2*b3 <= L - D + M*(1 - z12)
prob += p1 + a2 + b3 + a3 <= L - D + M*z13
prob += 1 - a2 - 2*b3 - 2*c3 <= L - D + M*(1 - z13)

# auxiliary variables for Fourier bound: A=max(a2+a4+a6,a3+a6,a5)
prob += a2 + a4 + a6 <= A
prob += a3 + a6 <= A
prob += a5 <= A
prob += b2 + b4 + b6 <= B
prob += b3 + b6 <= B
prob += b5 <= B
prob += c2 + c4 + c6 <= C
prob += c3 + c6 <= C
prob += c5 <= C

# Fourier bound
prob += A + B + C <= 4*L - 6*D

# Solve
prob.solve(pulp.PULP_CBC_CMD(msg=False))

# Output
print("Status:", pulp.LpStatus[prob.status])
for v in [a1,b1,c1,a2,b2,c2,a3,b3,c3,a4,b4,c4,a5,b5,c5,a6,b6,c6,D]:
    print(f"{v.name} = {v.value()}")

print("Maximal value of D:", pulp.value(prob.objective))
\end{lstlisting}

\end{document}